\documentclass
{article}
\usepackage{amsmath,amssymb,graphicx,theorem,latexsym}
\usepackage[all]{xy}
\xyoption{2cell}
\UseAllTwocells
\usepackage[bbgreekl]{mathbbol}

\DeclareFontFamily{U}{rsf}{}
\DeclareFontShape{U}{rsf}{m}{n}{
  <5> <6> rsfs5 <7> <8> <9> rsfs7 <10->  rsfs10}{}
\DeclareMathAlphabet{\mathscr}{U}{rsf}{m}{n}
\newcommand{\mycal}[1]{\mathscr{#1}}
\newtheorem{thm}{Theorem}[section]
\newtheorem{cor}[thm]{Corollary}
\newtheorem{lem}[thm]{Lemma}

\newtheorem{defn}[thm]{Definition}

\newtheorem{observation}[thm]{Observation}
\newtheorem{rem}[thm]{Remark}
{\theorembodyfont{\rmfamily} }
{\theorembodyfont{\rmfamily} }
\numberwithin{equation}{section}
\newcommand{\ignore}[1]{}

\newcommand{\delq}[1]{\frac{\partial}{\partial q_{#1}}}
\newcommand{\delp}[1]{\frac{\partial}{\partial p_{#1}}}

\newcommand{\Ch}{{\mathbb C}[[\hbar]]}
\newcommand{\bpi}{\boldsymbol{\Pi}}

\newcommand{\Ah}[1]{\mycal{A}_{#1}}
\newcommand{\bbD}{{\mathbb D}}

\newcommand{\bbX}{\mathbb{X}}

\newcommand{\cntrct}                
{\hspace{2pt}\raisebox{1pt}{\text{$\lrcorner$}}\hspace{2pt}}
\newcommand{\cntrctother}
{\hspace{2pt}\raisebox{1pt}{\text{$\llcorner$}}\hspace{2pt}}

\begin{document}
\title{Cohomology Groups of Deformations of Line Bundles on Complex Tori}
\author{O. Ben-Bassat \and  N. Solomon}

\date{August 28, 2009}

\maketitle

\begin{abstract}
The cohomology groups of line bundles over complex tori (or abelian varieties) are classically studied invariants of these spaces.  In this article, we compute the cohomology groups of line bundles over various holomorphic, non-commutative deformations of complex tori. Our analysis interpolates between two extreme cases.  The first case is a calculation of the space of (cohomological) theta functions for line bundles over constant, commutative deformations.  The second case is a calculation of the cohomologies of non-commutative deformations of degree-zero line bundles.
\end{abstract}


\tableofcontents

\section{Introduction}

Our goal in this article is to complete an analysis of the cohomology groups of line bundles on a non-commutative deformations of complex tori which began in \cite{bbp}.  In particular, that work studied a non-commutative, formal deformation of a complex torus $X$.   The sheaf of holomorphic functions on $X$ has a formal deformation $\Ah{\bpi}$ over the formal disk  in the direction of a holomorphic Poisson structure $\bpi$, where the multiplication law of $\Ah{\bpi}$ takes on the familiar Moyal form.   This is a particular example of a deformation quantization \cite{flato}, \cite{flato2},  \cite{Kontsevich}, \cite{kontsevich-course}, \cite{kontsevich-poisson}, \cite{NestTsygan}, \cite{moyal}.  In \cite{bbp}, the derived category of sheaves of $\Ah{\bpi}$-modules on a complex torus, $D^{b}(\Ah{\bpi}-mod)$, was shown to be equivalent via a Fourier-Mukai transform \cite{mukai} to the derived category of sheaves of modules over a formal gerbe on the dual torus $X^{\vee}$.

The {\it line bundles} referred to above are defined to be locally free, rank one $\Ah{\bpi}$-modules.
In this paper, we study the full subcategory of the $\mathbb{C}[[\hbar]]$ linear category $D^{b}(\Ah{\bpi}-mod)= D^{b}(\bbX_{\bpi})$ consisting of line bundles.
Therefore the goal of this paper is to compute the cohomology groups (which are $\mathbb{C}[[\hbar]]$-modules) of every line bundle $\mycal{L}$ over $\Ah{\bpi}$.
    Our motivation comes from two sources.  The first is an extension of Zharkov's \cite{Zharkov} work on cohomological theta functions (theta forms).  We expect this to play a role in checking the quantum background independence of certain deformations of the B-model.  In particular, we would like to extend Witten's work \cite{Witten} on background independence in string theory to the non-commutative holomorphic setting (a formal deformation of the B model)  featured for instance in \cite{Kapustin} \cite{block-dg} and \cite{bbp}.  The second source of motivation is a quest to understand in more detail the categories implicit in \cite{bbp}.  Specifically, we hope that after passing to algebroid stacks of quantizations, one can derive analytic descriptions \cite{Palamodov}, \cite{NestTsygan} of algebraic deformation quantizations \cite{Yekutieli1},\cite{Yekutieli2},  \cite{Yekutieli3}, \cite{Kontsevich},  \cite{vdb}, \cite{PoSc},  or non-commutative algebraic geometry on abelian varieties \cite{manin1}, \cite{manin2}.  In the future, we would also like to compare to work on deformed vector bundles in closely related categories as appear in  \cite{kaj1}, \cite{kaj2}, \cite{pol-sch}, \cite{pol1}, \cite{pol2}, \cite{calaque-halbout}, \cite{block-dg}, \cite{block-dg-next} and  \cite{block-daenzer}.  In \cite{KaSc} appears the first steps in a general study of modules over deformation quantizations on complex manifolds.
    Our main result is Theorem \ref{main} which expresses the $j$-th cohomology group of a $\Ah{\bpi}$ line bundle $\mycal{L}$ on the non-commutative complex torus $\bbX_{\bpi}= (X, \Ah{X, \bpi})$ in terms of (1) the $(j-1)$-st cohomology group of $L$, the classical line bundle given by the reduction modulo $\hbar$ of $\mycal{L}$ (2) A natural number $t^{0}$ determined (see Definition \ref{zerop}) by $\mycal{L}$ and (3) A corresponding element $l_{t^{0}} \in  \overline{V}^{\vee} = H^{1}(X,\mathcal{O})$ also determined by $\mycal{L}$.  The description of the series \[\hbar l_{1} + \hbar^{2}l_{2} + \hbar^{3}l_{3} + \cdots \] can be found in subsection \ref{LBaTD}.
\subsection{Acknowledgments}
We would like to thank Amnon Yekutieli for his helpful comments.

\subsection*{Notation and terminology}

\begin{description}
\item[$\Ah{S}$] a sheaf of associative flat $\Ch$ algebras on a complex manifold $S$
  satisfying $\Ah{S}/\hbar \cong \mathcal{O}_{S}$.
\item[$\Ch$] the complete local algebra of formal power series in
  $\hbar$.
\item[$\bbD$] the one dimensional formal disk.
\item[$V$] a complex vector space of dimension $g$.
\item[$\Lambda \subset V$] a free abelian subgroup of rank $2g$.
\item[$\bpi$] a holomorphic Poisson structure of constant rank $n$ on a complex manifold.
\item[$X$] a complex torus of dimension $g$, $X = V/\Lambda$.
\item[$V_{H,0}$] the complex subspace of $V$ annihilated by $H$.
\item[$g_{0}$] the dimension of $V_{0}$.
\item[$X_{H,0}$] the complex subtorus of $X$ annihilated by $H$, $X_{H,0} = V_{H,0}/(\Lambda \cap V_{H,0})$.
\item[$\iota$] the inclusion $V_{H,0} \to V$.
\item[$s$] a splitting of $\iota$.
\item[$\rho$] the quotient map $\rho:X \to X/X_{H,0}$.
\item[$p$] the quotient map $p:V \to V/\Lambda = X$.
\item[$(H,\chi)$] Appell-Humbert data defining a line bundle $L$, see a book such as \cite{lange}, \cite{mumford}, or \cite{polishchuk}.
\item[$\bbX_{\bpi}$] the Moyal quantization of the Poisson torus
  $(X,\bpi)$.
\item[$\mycal{L} = \mycal{L}_{((H,\chi);l(\hbar))}$] a line bundle on $\bbX_{\bpi}$, see equation \ref{eq:DefOfL}.
\item[$L$] a line bundle on $X$ given by the reduction of $\mycal{L}$ modulo $\hbar$.
\item[$\Phi$] the cocycle defining $\mycal{L}$, see equation \ref{eq:qAH}.
\item[$\phi$] the cocycle defining $\mycal{L}/\hbar^{t^{0}}$, see equation \ref{regularphi}.
\item[$\varphi$] the cocycle defining $L$.
\item[$l$] an element of $\overline{V}^{\vee}$.
\item[$l^{0}$] the image of $l$ under the projection to $\overline{V_{H,0}}^{\vee}$.
\item[$l(\hbar)$] an element of $\overline{V}^{\vee}[[\hbar]]$ given by \[\hbar l_{1} + \hbar^{2}l_{2}+ \hbar^{3}l_{3} + \cdots .\]
\item[$l(\hbar)^{0}$] an element of $\overline{V_{H,0}}^{\vee}[[\hbar]]$ given by \[\hbar l^{0}_{1} + \hbar^{2}l^{0}_{2}+ \hbar^{3}l^{0}_{3} + \cdots .\]
\item[$t$] a natural number depending on $\mycal{L}$, see definition \ref{tDef}.
\item[$t^{0}$] a natural number depending on $\mycal{L}$, see definition \ref{zerop}.
\item[$\overline{L}$] a line bundle on $X/X_{H,0}$, see Lemma \ref{lem:facts}.
\item[$\overline{\varphi}$] the cocycle defining $\overline{L}$
\item[$k$] an integer depending on $L$, see Lemma \ref{lem:facts}.
\end{description}

\section{Some Background} \label{sec:background}
First, let us recall the definition of complex tori and the non-commutative sheaf of algebras that we
will be using. For a more detailed discussion of the properties of
complex tori the reader may consult
\cite{mumford,lange,polishchuk}.  The description of the non-commutative sheaf of algebras is taken directly from \cite{bbp}.

A complex torus is a compact complex manifold $X$ which is isomorphic
to a quotient $V/\Lambda$, where $V$ is a $g$-dimensional complex
vector space and $\Lambda \subset V$ is a free abelian subgroup of
rank $2g$. Note that by construction $X$ has a natural structure of an
analytic group induced from the addition law on the vector space $V$.

Given a holomorphic Poisson structure \[
\bpi \in H^{0}(X, \wedge^{2}T_{X}) = \wedge^{2}V,
\]
 our non-commutative sheaf of algebras on $X$ is given by $\mathcal{O}_{X}[[\hbar]]$ as a sheaf of $\mathbb{C}[[\hbar]]-$modules along with a {\it Moyal product} $\star$ \cite{moyal,flato2} which we now describe.  We start by describing the
  standard Moyal product on $\mathcal{O}_{V}[[\hbar]]$ over the complex
  vector space $V$ equipped with $\bpi$ now viewed as a constant Poisson structure on $V$.
  \[\bpi \in \wedge^{2}V \subset H^{0}(V, \wedge^{2}T_{V}).
  \]
By the constancy, there are complex coordinates
  \[(q_1,\ldots, q_n, p_1,\ldots, p_n, c_1,\ldots, c_l)
  \] on
  $V$ so that the Poisson structure is diagonal, that is
\[
\bpi=\sum_{i=1}^{n} \delq{i}\wedge \delp{i}.
\]

With this notation we can now use $\bpi$ to define the
bidifferential operator $P$ by
\begin{equation} \label{eq:bidifferential}
P=\sum_i \left( \overleftarrow{\delq{i}}\,
\overrightarrow{\delp{i}} - \overleftarrow{\delp{i}}\,
\overrightarrow{\delq{i}}\right)
\end{equation}

Consider the sheaf $\mathcal{O}_V[[\hbar]]$ on $V$.
For any open set $U\subset V$, and any $f, g\in
\mathcal{O}(U)[[\hbar]]$ we define their Moyal product by
\[
f\star g=\sum_k\frac{\hbar^k}{k!} f\cdot P^k \cdot g=f\cdot\exp(\hbar
P)\cdot g = fg + \hbar\{ f, g \} + \cdots.
\]
Since the $\star$-product is defined by holomorphic bidifferential
operators it maps holomorphic functions to holomorphic functions.
Moreover, since bidifferential operators are local, the product
sheafifies. We denote the resulting sheaf of $\Ch$-algebras on $V$
by $\Ah{V,\bpi}$.
To define the Moyal quantization $(X,\Ah{X,\bpi})$ of a holomorphic
Poisson torus $(X,\bpi)$ we use the realization of $X$ as a
quotient $X = V/\Lambda$.  Let $p:V\to X$ be the covering
projection. Define the sheaf $\Ah{X,\bpi}$ of $\Ch$-algebras on $X$ as
follows. As a sheaf of $\mathbb{C}_{X}[[\hbar]]$-modules  it will
be just $\mathcal{O}_{X}[[\hbar]]$. To put a $\star$-product on
this sheaf one only has to use the natural identification
$\mathcal{O}_{X}[[\hbar]] := (p_{*}\mathcal{O}_{V}[[\hbar]])^{\Lambda}$
and note that the $\bpi$-Moyal product on $V$ is translation
invariant by construction. Explicitly the sections of $\Ah{X,\bpi}$
over $U\subset X$ can be described as the invariant sections
\begin{equation}\label{inv}
\Ah{X,\bpi}(U)=\Ah{V,\bpi}(p^{-1}(U))^\Lambda
\end{equation}
on the universal cover $V$. This is
well-defined since the Poisson structure $\bpi$ is constant and thus the
operator $P$ is translation invariant.

\subsection {A Review of Group Cohomology}
We now recall some basic definitions from group cohomology.  We describe the cohomology groups of various sheaves on the torus $V/\Lambda$ in terms of the group cohomology of $\Lambda$ acting on certain modules.
Let $G$ be a group and $M$ a $G$-module.  Denote the action of $G$ on $M$ by
\[(g,m) \mapsto g \cdot m.
\]
Recall \cite{mumford} that
the group cohomology differential
\[\delta: C^{p}(G, M)
\to C^{p+1}(G, M)
\]
is given by
\begin {equation}
\label {grcodif}
(\delta f)_{\lambda_{0}, \dots , \lambda_{p}} = \lambda_{0} \cdot (f_{\lambda_{1},\dots, \lambda_{p}}) + (\sum_{i=0}^{p-1} (-1)^{i+1} f_{\lambda_{0},\dots,\lambda_{i}+\lambda_{i+1}, \dots, \lambda_{p}}) + (-1)^{p+1}f_{\lambda_{0},\dots, \lambda_{p-1}}.
\end {equation}
Given a pairing
\[*: M \times N \to P
\]
of $G$-modules
the cup product on the level of cocycles is given by the map
\[C^{p}(G, M) \times C^{q}(G, N) \to C^{p+q}(G, P)
\]
\begin {equation}
\label {cup}
(f \cup g)_{\lambda_{0}, \dots , \lambda_{p+q-1}} = f_{\lambda_{0}, \dots, \lambda_{p-1}} * ((\lambda_{0}\dots \lambda_{p-1})g_{\lambda_{p+1}, \dots, \lambda_{p+q-1}})
.
\end {equation}
It is compatable with the differential and hence induces a cup product map
\[H^{p}(G, M) \otimes_{\mathbb{Z}}H^{q}(G, N) \overset{\cup}\to H^{p+q}(G,P).
\]
In the case where $M$, $N$, and $P$ are $R$-modules and the group action commutes with the action of $R$ and the pairing $*$ is $R$-bilinear, we actually get a cup product map
\[H^{p}(G, M) \otimes_{R} H^{q}(G, N) \overset{\cup} \to H^{p+q}(G,P).
\]
\begin{rem}\label{CaLe}
Given a sheaf of groups $\mathcal{S}$ on the torus $X$ with the property that $H^{j}(V,p^{-1}\mathcal{S})=0$ for $j>0$ we will often use the identification
\begin{equation}\label{important}H^{i}(X,\mathcal{S}) \cong H^{i}(\Lambda, p^{-1}\mathcal{S}(V))
.\end{equation}  This isomorphism comes from the collapse of the Cartan-Leray spectral sequence.
\end{rem}
The isomorphisms in (\ref{important}) are natural in $\mathcal{S}$ and the cup product we have explained above gives the cup product in sheaf cohomology.  The action of $\Lambda$ on $p^{-1}\mathcal{S}(V)$ is by translation.
In this paper, we work with $\mathbb{C}[[\hbar]]-$modules for the group $\Lambda$ such that the action of $\Lambda$ commutes with the multiplication by elements of $\mathbb{C}[[\hbar]]$.  Therefore the cohomology groups are in a natural way $\mathbb{C}[[\hbar]]-$modules and to compute them it will often be useful to use the following observation.
\begin{observation}\label{structure}  The structure theorem for finitely generated modules over the principal ideal domain $\mathbb{C}[[\hbar]]$ shows that every finitely
generated $\mathbb{C}[[\hbar]]-$module is determined uniquely up to isomorphism by its associated graded vector space.
\end{observation}
\subsection{Line Bundles and Their Deformations}\label{LBaTD}
Recall from \cite{bbp} that a {\it line bundle} on $\bbX_{\bpi}$ is a locally free, rank one left $\Ah{\bpi}$-module in the classical topology on $X$.  In \cite{bbp}, a type of Appell-Humbert theorem was proven, whereby all line bundles on $\bbX_{\bpi}$ were classified and constructed.  Namely, a correspondence was established between equivalence classes of line bundles on $\bbX_{\bpi}$ and Appell-Humbert data consisting of certain triples $((H,\chi);l(\hbar))$.  Here, $H$ is an element of the Neron Severi group of $X$, $\chi:\Lambda \to U(1)$ is a semi-character for $H$ and $l(\hbar) \in \hbar(\overline{V}^{\vee}[[\hbar]])$.

The pair $(H,\chi)$ is the classical Appell-Humbert data corresponding to equivalence classes of classical line bundles $L$, while $l(\hbar)$ describes $\mycal{L}$ as an iterated extension of $L = \mycal{L}/\hbar$ by itself.  Line bundles $\mycal{L}$ on $\bbX_{\bpi}$, are classified \cite{bbp} up to equivalence by triples $((H,\chi);l(\hbar))$  satisfying the equation
\begin{equation}\label{condition}
H \cntrct
\bpi \cntrctother H  = 0.
\end{equation}

   In order to explain this condition notice that since we are working over a torus, we can consider $H \in (V^{\vee} \otimes \overline{V}^{\vee}) \cap Alt^{2}(\Lambda, \mathbb{Z})$.  In other words, $H$ is a Hermetian form on $V$ which satisfies $\text{Im} H(\Lambda, \Lambda) \in \mathbb{Z}$.  The contraction in equation (\ref{condition}) can be described as the usual contraction of $H \wedge H \in (\wedge^{2} V^{\vee})\otimes (\wedge^{2}\overline{V}^{\vee})$ with $\bpi$.

We now describe how line bundles on $\bbX_{\bpi}$ can be constructed from Appell-Humbert data.  Equivalence classes of line bundles on $\bbX_{\bpi}$ are in one to one correspondence with elements of the pointed set \[H^{1}(X,\Ah{X,\bpi}^{\times}) \cong H^{1}(\Lambda, \Ah{V,\bpi}^{\times}(V)).\]  The existence of the above isomorphism follows immediately from Remark \ref{CaLe} and the exponential short exact sequence of sheaves of groups on $V$:
\[0 \to \mathbb{Z} \to \mathcal{O}_{V}[[\hbar]] \to \Ah{V, \bpi}^{\times} \to 1.
\]
  We will denote by $T_{\lambda}$ the automorphism of $V$ given by $v \mapsto v+ \lambda$.
The set $Z^{1}(\Lambda, \Ah{\bpi}^{\times}(V))$ consists of maps $\Phi:\Lambda \to \Ah{V,\bpi}^{\times}(V)$ satisfying
\[\Phi_{\lambda_{2}} \star ( \Phi_{\lambda_{1}} \circ T_{\lambda_{2}}) = \Phi_{\lambda_{1}+\lambda_{2}}.
\]
For each triple $((H,\chi);l(\hbar))$ satisfying equation (\ref{condition}) we define a line bundle $\mycal{L}_{((H,\chi);l(\hbar))}$.
First we construct an element $\Phi = \Phi_{((H,\chi),\, l(\hbar))}$ in $Z^{1}(\Lambda, \Ah{\bpi}^{\times}(V))$ given by
\begin{equation} \label{eq:qAH}
\Phi_{((H,\chi),\, l(\hbar))}(\lambda)(v) = \chi(\lambda)
\exp\left(\pi H(v, \lambda) +
\frac{\pi}{2}H(\lambda, \lambda) +
\sum_{m=1}^{\infty} \hbar^{j}\pi \langle l_{m},\lambda \rangle \right).
\end{equation}

This is a non-commutative deformation of the factor of automorphy $\varphi_{(H,\chi)}$ corresponding to the line bundle $L = L_{(H,\chi)}$,

\begin{equation} \label{eq:qAHcl}
\varphi_{(H,\chi)}(\lambda)(v) = \chi(\lambda)
\exp\left(\pi H(v, \lambda) +
\frac{\pi}{2}H(\lambda, \lambda)\right).
\end{equation}
That is, we have \[\varphi = \Phi \mod \hbar.\]
Finally, we define
\begin{equation} \label{eq:DefOfL}
\mycal{L} = \mycal{L}_{((H,\chi);l(\hbar))} \subset p_{*} \mathcal{O}_{V}[[\hbar]],
\end{equation}
to be the subsheaf of $p_{*} \mathcal{O}_{V}[[\hbar]]$ consisting for small enough $U \subset X$ of
\[\{f \in \mathcal{O}_{V}(p^{-1}(U))[[\hbar]] | f \circ t_{\lambda} = f \star \Phi_{\lambda} \quad \forall \lambda \in \Lambda \}
.\]
Here $t_{\lambda}$ is the automorphism of $p^{-1}(U)$ given by $w \mapsto w+ \lambda$.
\subsection{The Twisted Action}
We wish to compute the $\mathbb{C}[[\hbar]]-$module $H^{j}(X,\mycal{L})$.  In the Appendix, we use the canonical isomorphism $p^{-1}\mycal{L} \cong \mathcal{O}_{V}[[\hbar]]$ to conjugate the translation action $A$ of $\Lambda$ on $p^{-1}\mycal{L}(V)$ to obtain an action $A^{\Phi}$ of $\Lambda$ on $\mathcal{O}(V)[[\hbar]]$.  The formula derived in the Appendix is
\begin{equation}A^{\Phi}_{\lambda}(g) = (A_{\lambda}(g)) \star \Phi_{\lambda}^{-1}
.\end{equation}
In other words, the action $A^{\Phi}_{\lambda}$ on $H^{0}(V, \mathcal{O}_{V} [[\hbar]])$ is given by
\begin {equation}
\label {trac}
g(v) \mapsto g(v +\lambda) \star \Phi_{\lambda}(v)^{-1}.
\end {equation}

 We denote the corresponding cohomology groups corresponding to the $A^{\Phi}$ action by $H^{i}(\Lambda,\mathcal{O}(V)[[\hbar]], \Phi)$. By the collapse of the Cartan-Leray spectral sequence \ref{important} (which happens because $H^{j}(V, p^{-1}\mycal{L}) = H^{j}(V, \mathcal{O}_{V}[[\hbar]])=0$ for $j>0$) we have isomorphisms of $\mathbb{C}[[\hbar]]-$modules
\begin {equation}
\label {CarLer}
H^{i}(X,\mycal{L}) \cong H^{i}(\Lambda,p^{-1}\mycal{L}(V)) \cong H^{i}(\Lambda,\mathcal{O}(V)[[\hbar]], \Phi).
\end {equation}
Thus, it remains to compute the cohomology groups $H^{i}(\Lambda,\mathcal{O}(V)[[\hbar]], \Phi)$ .

\section{The Main Computation}

\subsection{The Spectral Sequence}
Consider the filtration on $\mycal{L}$ defined by
$F^n \mycal{L} = \hbar^n \mycal{L}$.  This filtration induces a filtration on the standard complex computing group cohomology.  Then by applying the {\it Spectral Sequence of a Filtration} of section 5.4 of \cite{Weibel} to this filtered complex, we get a spectral sequence
\begin{equation}\label{GradSpec}E_1^{p,q} := H^{p+q}(\Lambda,Gr_{p} (\mathcal{O}_{V}(V)[[\hbar]]), \varphi) \Rightarrow H^{p+q}(\Lambda,\mathcal{O}_{V}(V)[[\hbar]],\Phi)\end{equation}
or equivalently
$$E_1^{p,q} := H^{p+q}(X,Gr_{p} \mycal{L}) \Rightarrow H^{p+q}(X,\mycal{L}),$$
with differentials
\[d_{j}^{p,q}: E_{j}^{p,q} \to E_{j}^{p+j,q-j+1}.
\]
In order to do computations, we will use the group cohomology version but in order to simplify notation we will write the terms with sheaf cohomology.  It is important to note that the appearance of $\varphi$ in (\ref{GradSpec}) is due to the fact for each $\lambda \in \Lambda$ that the $p-$th graded part of the action of $A^{\Phi}_{\lambda}$ is simply the action of $A^{\varphi}_{\lambda}$.

Notice that \[Gr_p \mycal{L} = F^p\mycal{L}/F^{p+1}\mycal{L} \cong L,\] where $L$ is a (classical) line bundle on $X$.  Therefore the $E_{1}$ term looks like
\def\csg#1{\save[].[ddddddrrrrr]!C*+<2pc>[F-,]\frm{}\restore}
\[
\xymatrix@=1.5pc{
q=3 & \csg1 \ \ \ \ \ \ \ \vdots & \vdots & \vdots & \vdots & \vdots &\\
q=2 &
H^{2}(X, L)   &
H^{3}(X,L) &
H^{4}(X,L) & H^{5}(X,L) &
H^{6}(X,L) & \cdots \\
q=1 &
H^{1}(X, L)   &
H^{2}(X,L) &
H^{3}(X,L) & H^{4}(X,L) &
H^{5}(X,L) & \cdots \\
q= 0 & H^{0}(X,L)   &
H^{1}(X,L) &
H^{2}(X,L) & H^{3}(X,L) &
H^{4}(X,L) & \cdots \\
q=-1 &    &
H^{0}(X,L)  &
H^{1}(X,L) & H^{2}(X,L) & H^{3}(X,L) & \cdots \\
q= -2 &
&   &
H^{0}(X,L) & H^{1}(X,L) & H^{2}(X,L) & \cdots \\
q=-3 & & & & \ddots & \ddots & \ddots \\
 & p= 0 & p= 1 & p= 2 & p=3 & p=4
}
\]

\begin{defn}\label{DefOfSpectralSequence}
The spectral sequence associated to $X$ and $\mycal{L}$ as above will be denoted $(E(X,\mycal{L}),d(X,\mycal{L}))$ or abbreviated by $(E,d)$ when $X$ and  $\mycal{L}$ are clear.
\end{defn}
\begin{thm}
\label {specseqthm}
When the spectral sequence $(E,d)$ converges, there exist canonical isomorphisms
\[Gr^{p}(H^{q}(X,\mycal{L})) \cong E_{\infty}^{p,p+q},
\]
for given integers $p,q$.
\end{thm}

{\bf Proof.}
This is a Corollary of Theorem 5.5.10 of \cite{Weibel}.

\ \hfill $\Box$

\begin {defn}\label{SpecShift}
Let $a,b$ be integers. For any spectral sequence $(E,d)$, we define the shifted spectral sequence
$(E,d)[a,b]$ by
$$(E[a,b]^{p,q}, d[a,b]^{p,q}):= (E^{p-a,q-b}, d^{p-a,q-b}).$$
\end {defn}

\begin {defn}\label{tDef}
Let $\mycal{L} = \mycal{L}_{((0,1);l(\hbar))}$ be a line bundle on $\bbX_{\bpi}$ such that $\mycal{L}/\hbar$ is trivial.
 we define
$$t=t_{\mycal L} := \min \{s, l_{s} \ne 0\},$$
or $\infty$ if the minimum is not obtained. \\
\end {defn}
\begin {lem} \label{someLemma1}
We have
\[
d_{j}(\bullet) =  -\pi l_{j} \cup \bullet \ \ \ \ \ \ \text{for} \ \ \ \ \ \ j \leq t .
\]
\end {lem}
{\bf Proof.}
    The proof is via induction on $j$.
 The base case of the induction is the case $j=1$.  Both the base case and the induction step follow from the following.
    \ \hfill $\Box$

 \begin{lem} \label{someLemma2} If $d_{s}=0$ and $l_{s}=0$ for all $s<j$ then $d_{j}(\bullet) =  -\pi l_{j} \cup \bullet $.
 \end{lem}
 {\bf Proof.}
    As $d_s = 0$ for $s<j$ we deduce $E_j^{p,q} = E_1^{p,q} = H^{p+q}(X,Gr_{p} \mycal{L}) = H^{p+q}(X,L)$ and so the differential $d_{j}$ is a map \[d_{j}: H^{p+q}(X,Gr_{p} \mycal{L}) \to H^{p+q+1}(X,Gr_{p+j} \mycal{L}).
    \]

    Given $[\xi] \in E_j^{p,q}  \cong H^{p+q}(\Lambda,\mathcal{O}(V),\varphi)$, we compute $d_j^{p,q}(\xi)$ according to the prescription in Appendix 2, item 4. Let $i=p+q$.  Consider the element $\xi$ of $Z^{i}(\Lambda,\mathcal{O}(V),\varphi)$.  The boundary under the map $\delta^{\Phi}$ then lies in $C^{i+1}(\Lambda, \hbar^{j}\mycal{A}(V),\Phi).$  Considering this element modulo $\hbar^{j+1}$ finally gives the desired class in $Z^{i+1}(\Lambda, Gr_{j}\mycal{A}(V),\varphi)$.\\
    Let $\xi \in Z^{i}(\Lambda,\mathcal{O}(V),\varphi)$.

    \begin {equation}
    \begin {array}{ll}
    \label {form}
    (\delta^{\Phi} \xi)_{\lambda_{0}, \dots , \lambda_{i}} = A_{\lambda_{0}}^\Phi (\xi_{\lambda_{1},\dots, \lambda_{i}}) + (\sum_{c=0}^{i-1} (-1)^{c+1} \xi_{\lambda_{0},\dots,\lambda_{c}+\lambda_{c+1}, \dots, \lambda_{i}}) + (-1)^{i+1}\xi_{\lambda_{0},\dots, \lambda_{i-1}}  \cr\cr
    = (A_{\lambda_{0}}^\Phi -A^{\varphi}_{\lambda_{0}}) (\xi_{\lambda_{1},\dots, \lambda_{i}}) + A^{\varphi}_{\lambda_{0}} (\xi_{\lambda_{1},\dots, \lambda_{i}}) + (\sum_{c=0}^{i-1} (-1)^{c+1} \xi_{\lambda_{0},\dots,\lambda_{c}+\lambda_{c+1}, \dots, \lambda_{i}}) + (-1)^{i+1}\xi_{\lambda_{0},\dots, \lambda_{i-1}}

    \cr\cr = (\xi_{\lambda_{1},\dots, \lambda_{i}} \circ T_{\lambda_{0}}) \star {(\Phi_{\lambda_{0}}^{-1} - \varphi_{\lambda_{0}}^{-1})} + (\delta^{\varphi} \xi)_{\lambda_{0}, \dots , \lambda_{i}}
    \cr\cr = (\xi_{\lambda_{1},\dots, \lambda_{i}} \circ T_{\lambda_{0}}) \star {(\Phi_{\lambda_{0}}^{-1} - \varphi_{\lambda_{0}}^{-1})}

    .
    \end {array}
    \end {equation}
    Now, we divide the resulting element by $\hbar^j$.
    and take the resulting element modulo $\hbar$.
    We have
    $$(\Phi_{\lambda_{0}}^{-1} - \varphi_{\lambda_{0}}^{-1}) =  \varphi_{\lambda_{0}}^{-1}\left(\exp\left(-\sum_{m=1}^{\infty} \hbar^{m}\pi \langle l_{m},\lambda_0 \rangle \right)-1 \right).$$
    Thus, modulo $\hbar$ we have
    \begin {equation}
    \begin {array}{ll}
    \frac {(\Phi_{\lambda_{0}}^{-1} - \varphi_{\lambda_{0}}^{-1})} {\hbar^j} \equiv  -\varphi_{\lambda_{0}}^{-1} \cdot \pi \langle l_{j},\lambda_0 \rangle \mod \hbar.
    \end {array}
    \end {equation}
    The resulting element is therefore
    \begin {equation}
    \begin {array}{ll}
    -\varphi_{\lambda_{0}}^{-1} \cdot \pi \langle l_{j},\lambda_0 \rangle ( \xi_{\lambda_{1},\dots, \lambda_{i}}\circ T_{\lambda_{0}}) = -\pi \langle l_{j},\lambda_0 \rangle A^{\varphi}_{\lambda_{0}}(\xi_{\lambda_{1},\dots, \lambda_{i}}) .
    \end {array}
    \end {equation}
    By reference to formula (\ref{cup}) this element in $Z^{i+1}(\Lambda,\mathcal{O}(V),\varphi)$ represents the image of the cup product
    \[H^{1}(X,\mathcal{O}) \otimes H^{i}(X,L) \to H^{i+1}(X,L).
    \]
    evaluated on $[\lambda \mapsto -\pi \langle l_{j},\lambda \rangle]$ and $[\xi]$.  In other words
     \[d_{t}([\xi]) = -\pi l_{j} \cup [\xi]
    \]
    as claimed.
    \ \hfill $\Box$

   The preceding analysis shows that the spectral sequence (\ref{GradSpec}) satisfies $E_{t} = E_{1}$.

\subsection{The Case of a Deformation of the Trivial Line Bundle}
Let $\mycal{L}$ be a line bundle on $\bbX_{\bpi}$.   We now compute the term $E_{t+1}$ of the spectral sequence (\ref{GradSpec}) in the case that $\mycal{L}/\hbar = L \cong \mathcal{O}$.  Note that we have $\mycal{L} = \mycal{L}_{((0,1,);l(\hbar))}$.
Recall that the canonical isomorphisms $H^{a}(X,\mathcal{O}) \cong \wedge^{a} \overline{V}^{\vee}$ fit into a commutative diagram \cite{mumford}
\[
\xymatrix@C+2pc{
H^{a}(X,\mathcal{O}) \otimes  H^{b}(X,\mathcal{O}) \ar[d]
\ar[r]_-{\cup}
 &  H^{a+b}(X,\mathcal{O}) \ar[d]
  \\
\wedge^{a} \overline{V}^{\vee} \otimes \wedge^{b} \overline{V}^{\vee}
\ar[r]_-{\wedge}
&   \wedge^{a+b} \overline{V}^{\vee}.
}
\]
This implies that the map $d_{t}$ on $E_{t}$ becomes several copies of a truncated version of the Koszul sequence $(\wedge^{\bullet} \overline{V}^{\vee}, (\bullet)\wedge -\pi l_{t})$ for the wedge product by the element $-\pi l_{t}$.  These copies look like
\[(E_{t}^{p,q} = \wedge^{p+q} \overline{V}^{\vee}) \to (E_{t}^{p+t,q-t+1} = \wedge^{p+q+1} \overline{V}^{\vee}) \to (E_{t}^{p+2t,q-2t+2} = \wedge^{p+q+2} \overline{V}^{\vee}) \to \cdots \to \wedge^{g} \overline{V}^{\vee} \to 0
\]
where $0 \leq p <t$ and $-p \leq q \leq g-p$.
The (untruncated) Koszul sequence $(\wedge^{\bullet} \overline{V}^{\vee}, (\bullet)\wedge -\pi l_{t})$ is exact.  Therefore, the cohomology of the  truncated Koszul sequence is just the kernel of the wedge product with $l_{t}$.   Therefore \[E_{t+1}^{p,q} =0 \ \ \text{       if        } \ \
p \ge t  \] and
\[E_{t+1}^{p,q} =\ker (\bullet \wedge l_{t}:\wedge^{p+q} \overline{V}^{\vee} \to \wedge^{p+q+1} \overline{V}^{\vee}) = \text{im} (\bullet \wedge l_{t}:\wedge^{p+q-1} \overline{V}^{\vee} \to \wedge^{p+q} \overline{V}^{\vee}) = \wedge^{p+q}\overline{V}^{\vee}/<l_{t}>
\ \  \text{     if     } \ \ p < t.  \]
It is clear then that $d_{t+1}$ is identically 0, and therefore the spectral sequence degenerates at $E_{t+1}$.  Therefore the $E_{\infty}$ term looks as follows:

\def\csg#1{\save[].[ddddddddrrrrr]!C*+<3.5pc,1pc>[F-,]\frm{}\restore}
\[
\xymatrix@=.5pc{
q=3 & \csg1 \ \ \ \ \ \ \vdots & \vdots & \vdots & \vdots & \vdots &\\
q=2 &
\wedge^{2}\overline{V}^{\vee}/<l_{t}>   &
\wedge^{3}\overline{V}^{\vee}/<l_{t}> &
\dots & \wedge^{t+1}\overline{V}^{\vee}/<l_{t}>  &
0 & \cdots \\
q=1 &
\overline{V}^{\vee}/<l_{t}>   &
\wedge^{2}\overline{V}^{\vee}/<l_{t}> &
\dots & \wedge^{t}\overline{V}^{\vee}/<l_{t}> &
0 & \cdots \\
q= 0 & 0  &
\overline{V}^{\vee}/<l_{t}>  &
\dots & \wedge^{t-1}\overline{V}^{\vee}/<l_{t}>  &
0 & \cdots \\
q=-1 &    &
0  &
\ddots & \vdots  & 0 & \cdots \\
\vdots &
&   &
\ddots & \overline{V}^{\vee}/<l_{t}> & 0 & \cdots \\
q=-t+1 & & & &0& 0 & \ddots \\
q=-t & & & &  & 0 & \ddots \\
 \vdots & & & &  &  & \ddots \\
 & p= 0 & p= 1 & p= 2 & p=t-1 & p=t & \cdots
}
\]
So we have \[E_{\infty}^{p,q}= \ker(\bullet \cup l_{t}: H^{p+q}(\Lambda,\mathcal{O}(V)) \to H^{p+q+1}(\Lambda,\mathcal{O}(V)))\] for $p<t$ and $0$ otherwise.

As the spectral sequence converges, Theorem \ref{specseqthm}, together with Observation \ref{structure} implies

\begin{lem}\label{cupl}Let $\mycal{L} = \mycal{L}_{((0,1,);l(\hbar))}$ be a line bundle on $\bbX_{\bpi}$ such that $\mycal{L}/\hbar \cong \mathcal{O}$ and $l(\hbar) \neq 0$.  Then there are isomorphisms of $\mathbb{C}[[\hbar]]-$modules
\[H^{j}(\Lambda,\mathcal{O}(V),\Phi) \cong \mathbb{C}[\hbar]/(\hbar^{t}) \otimes_{\mathbb{C}}(l_{t} \cup H^{j-1}(\Lambda,\mathcal{O}(V)))
\]
or equivalently
\begin {equation}
\label {iso1}
H^{j}(X,\mycal{L}) \cong \mathbb{C}[\hbar]/(\hbar^{t}) \otimes_{\mathbb{C}}(l_{t} \cup H^{j-1}(X,\mathcal{O})),
\end {equation}
where the $\mathbb{C}[[\hbar]]$ structure is inherited from the $\mathbb{C}[\hbar]/(\hbar^{t})$ term and the convention is that $H^{-1}(\Lambda,\mathcal{O}(V)) = 0$.  The term $l_{t} \cup H^{p-1}(\Lambda,\mathcal{O}(V))$ is the image under the linear map
\[H^{j-1}(\Lambda,\mathcal{O}(V)) \to H^{j}(\Lambda,\mathcal{O}(V))
\]
given by taking the cup product with $l_{t}$.
\end{lem}

\ \hfill $\Box$

\noindent
Therefore
\begin{cor}
 In the situation of Lemma \ref{cupl} above we have \[dim_{\mathbb{C}} H^{j}(X,\mycal{L}) = t {{g-1}\choose{j-1}}.
\]
\end{cor}

\ \hfill $\Box$

\subsection{The General Case}
In order to analyze the differential in the spectral sequence in the general case, we have to carefully combine the two extreme cases.   Let $\mycal{L}=\mycal{L}_{((H,\chi);l(\hbar)}$ be a line bundle on $\bbX_{\bpi}$.  Let $\iota: X_{H,0} \subset X$ be the degeneracy locus of $H$.  It is a sub torus of $X$ and $\rho: X \to X/X_{H,0}$ is a principal $X_{H,0}$ bundle.  Let $g_{0}$ be the dimension of $X_{H,0}$.  Fix, once and for all, a splitting $s:V \to V_{H,0}$ of the inclusion $V_{H,0} \to V$.  We use the same letter to denote the corresponding splitting $\wedge^{j}V \to \wedge^{j}V_{H,0}$ and $s^{\vee}$ to denote its complex conjugated dual $s^{\vee}: \wedge^{j}\overline{V_{H,0}}^{\vee} \to \wedge^{j}\overline{V}^{\vee}$, which are equivalently thought of as maps
\begin {equation}
\label {sdual}
s^{\vee}: H^{j}(X_{H,0},\mathcal{O}) \to H^{j}(X,\mathcal{O}).
\end {equation}
We also have 
$\iota^*: H^{j}(X,\mathcal{O}) \to H^{j}(X_{H,0},\mathcal{O})$.  We also notice that we have a split short exact sequence
\begin {equation}
\label {split}
0 \to H^1(X/X_{H,0},\mathcal O) \overset{\rho^*}\to H^1(X,\mathcal O) \overset{\iota^*} \to H^{1}(X_{H,0},\mathcal{O})\to 0,
\end {equation}
where the splitting is given by $s^\vee$.

\begin {defn}
\label {spl}
For $a \in H^j(X,\mathcal O)$ we denote by $a^0$ its image under the morphism
$$H^j(X,\mathcal O) \overset{\iota^*} \to H^{j}(X_{H,0},\mathcal{O}).$$
\end {defn}

\begin {defn}
\label {zerop}
Given a line bundle $\mycal L $ on $\bbX_{\bpi}$ we define
$$t^0=t^0_{\mycal L} := \min \{s, l^{0}_{s} \ne 0\},$$
or $\infty$ if the minimum is not obtained. We also let $l(\hbar)^0 = \sum_{i=1}^\infty l_i^0 \hbar^i.$
\end {defn}
From \cite{lange} we learn that

\begin{lem}\label{lem:facts}
There exists a line bundle $\overline{L}$ on $X/X_{H,0}$ such that
\begin{enumerate}
\item
The restriction of $L$ to $X_{H,0}$ is a degree zero line bundle which can be considered a restriction from $X$ to $X_{H,0}$ of a degree zero line bundle $P$ on
$X$.
\item
There is an isomorphism $L \cong P \otimes \rho^{*} \overline{L}$.
\item
There is a unique integer $k$ such that  $H^{k}(X/X_{H,0},\overline{L}) \neq  0$.
\item
If $P \cong \mathcal O$, then the only non-zero cohomology groups of $L$ have dimensions, for $0 \leq i \leq g_{0}$:
\begin{equation}\label{eqn:DimResult} h^{i+k}(X,L) = h^{k}(X/X_{H,0},\overline{L}){{g_{0}}\choose{i}}.
\end{equation}
Otherwise, $h^j(X,L)=0$ for all $j$.
\item
If $P \cong \mathcal O$, then the pullback map
$$\rho^*: H^{k}(X/X_{H,0},\overline{L}) \to H^{k}(X,L)$$
is an isomorphism.
\end{enumerate}
\end{lem}

\ \hfill $\Box$

Equation 3.5.1(1) in \cite{lanbir} implies 

\begin {lem}
\label {somsur}
The cup product map 
\begin{equation*}
\cup: H^{k}(X, L) \otimes H^{i}(X,\mathcal{O}) \to H^{k+i}(X,L)
\end{equation*}
is surjective.
\end {lem}

\begin {lem}
\label {somiso}
Let $\xi \in H^{k}(X/X_{H,0},\overline{L}), \ a \in H^{i}(X, \mathcal{O})$. Then
$$\rho^*(\xi) \cup a = \rho^*(\xi) \cup s^\vee(a^0).$$
\end {lem} 

{\bf Proof.}
We notice that
\[a-s^{\vee}(a^{0}) \in \rho^{*}(H^{1}(X/X_{H,0},\mathcal{O}))\cup H^{i-1}(X, \mathcal{O}) \subset H^{i}(X, \mathcal{O}).
\]

On the other hand, the cup product map
\[H^{k}(X/X_{H,0}, \overline{L})\otimes H^{1}(X/X_{H,0},\mathcal{O}) \overset{\cup} \to H^{k+1}(X/X_{H,0}, \overline{L})
\]
vanishes by Lemma \ref{lem:facts}.3.
This completes the proof.
\ \hfill $\Box$
\begin {defn}
For a line bundle $L$ on $X$, the integer appearing in Lemma \ref{lem:facts}.3 will be denoted by $k=k_L$ (we omit the subscript $L$ if there is no ambiguity).
\end {defn}
\begin{lem} \cite{lange}
\label{green}
When $L$ restricts to a non trivial line bundle on $X_{H,0}$ then $H^{l}(X,L)$ vanishes for all $l$.  If the restriction of $L$ to $X_{H,0}$ is isomorphic to $\mathcal{O}_{X_{H,0}}$ (so $P\cong \mathcal{O}_{X}$), then $\rho^{*} \overline{L} = L$  (see Lemma \ref{lem:facts}) and the map
\begin{equation}\label{identify}
H^{k}(X/X_{H,0},\overline{L}) \otimes H^{i}(X_{H,0},\mathcal{O}) \to H^{k+i}(X,L)
\end{equation}
defined by
\[b \otimes a \mapsto \rho^{*}(b) \cup s^{\vee}(a)
\]
is an isomorphism.  Consider a set $\{b^r \}$ of elements of $Z^{k}(\Lambda/(V_{H,0} \cap \Lambda), \mathcal{O}(V/V_{H,0}), \overline{\varphi})$ whose cohomology classes are a basis for $H^{k}(X/X_{H,0},\overline{L})$. Let $a^{m_{1}},\dots, a^{m_{i}}$ be a basis for $\overline{V_{H,0}}^{\vee}$.
If we define
\[a^{I} = a^{m_{1}} \cup \cdots \cup a^{m_{i}}
\]
for $1 \leq m_{1} < \dots m_{i} \leq g_{0}$ then the following collection
\begin{equation}\label{finalExplicitNol}
\left(b^{r}_{\rho(\lambda_{1}), \dots , \rho(\lambda_{k})} \circ \rho \right) \left(s^{\vee}a^{I}_{\lambda_{k+1}, \dots , \lambda_{k+i}}\right)
\end{equation}
are $h^{k}(X,L){{g_{0}}\choose{i}}$ elements in $Z^{k+i}(\Lambda,\mathcal{O}(V),\varphi)$ whose cohomology classes form a basis for $H^{k+i}(X,L)$.
\end{lem}
\noindent{\bf Proof.}
Since the dimensions of both vector spaces are the same (by Lemma \ref{lem:facts}.4), it is enough to show that the morphism above is surjective.
By Lemma \ref{lem:facts}.5, it is enough to show that the morphism
\begin{equation}\label{identify}
H^{k}(X,L) \otimes H^{i}(X_{H,0},\mathcal{O}) \to H^{k+i}(X,L)
\end{equation}
defined by
\[b \otimes a \mapsto b \cup s^{\vee}(a)
\]
is surjective.
But this follows by Lemmas \ref{somsur} and \ref{somiso}.

\ \hfill $\Box$

\begin {lem}
\label{cupprodiso}
For $l \in H^{1}(X,\mathcal{O})$, $a \in
H^{i}(X_{H,0},\mathcal{O})$ and $b \in H^{k}(X/X_{H,0}, \overline{L})$,
we have
\[(\rho^{*}(b) \cup s^{\vee}(a)) \cup l =  \rho^{*}(b) \cup s^\vee(a \cup l^{0}).
\]
\end {lem}
{\bf Proof.}
Lemma \ref{somiso} implies that 

\begin{equation}\label{CompatOfLs}(\rho^{*}(b) \cup s^{\vee}(a)) \cup l =  \rho^{*}(b) \cup s^\vee(a) \cup s^\vee(l^{0}).
\end{equation}
As 
\begin {equation}
\label {assoc}
s^{\vee}(a \cup l^{0}) = s^{\vee}(a) \cup s^{\vee}(l^{0}),
\end {equation}
the Lemma follows.
\ \hfill $\Box$

Recall from Definition \ref{DefOfSpectralSequence} the spectral sequence $(E(X,\mycal{L}), d(X,\mycal{L}))$ associated to any line bundle $\mycal{L}$ on $\bbX_{\bpi}$, together with Definition \ref{SpecShift} of the shifted spectral sequence in \ref{SpecShift}.
\begin {defn}
Let\[\tilde{E}(X,\mycal{L}) = E(X_{H,0},\mycal{L}_{((0,1);l(\hbar)^{0})})[0,k] \otimes H^{k}(X/X_{H,0},\overline{L}).
\]
and
\[\tilde{d}(X,\mycal{L}) = d(X_{H,0},\mycal{L}_{((0,1);l(\hbar)^{0})})\otimes \text{id}.
\]
\end {defn}

As the tensor product over $\mathbb{C}$ is exact, this is a spectral sequence.
\begin{lem}  There is an isomorphism of spectral sequences

\[(E(X,\mycal{L}), d(X,\mycal{L})) \cong (\tilde{E}(X,\mycal{L}), \tilde{d}(X,\mycal{L})).
\]

\end{lem}

{\bf Proof.}
The identification \ref{identify} sets up an isomorphism
\begin{equation}\label{EIsom}E(X,\mycal{L}) \cong \tilde{E}(X,\mycal{L}).
\end{equation}
We now show that $d(X,\mycal{L})$ and  $\tilde{d}(X,\mycal{L})$ correspond to one another under this identification.
The proof is done in two steps.
\begin {itemize}
\item 
Step $1$: We show that $d_j$ corresponds to $\tilde d_j$ under \ref{EIsom} for $j \le t^0$.\\
The proof is via induction on $j$.
The case $j=1$ follows by Lemma \ref{cupprodiso}.
We assume the claim holds true for $J<j$ and prove the claim for $j$.
For $b<t^0$, Lemma \ref{someLemma2} implies that
$\tilde d_b$ is given by cupping with $-\pi l_{b}^0$, and thus vanishes.
By the induction hypothesis for $J<j$, we deduce that $d_J = 0$.
Again using Lemma \ref{someLemma2}, we deduce that $d_j(\bullet) = -\pi l_j \cup (\bullet)$.
Another application of Lemma \ref{cupprodiso} completes the induction step.
\item Step $2$: We show that $d_j$ corresponds to $\tilde d_j$ under \ref{EIsom} for $j > t^0$.\\
Since the $t_0+1$ sheets of both spectral sequences are concentrated in a block of width $t_0$, both $d_j$ and $\tilde d_j$ are 0 for $j>t_0$.
\end {itemize}

This completes the proof.
\ \hfill $\Box$
\begin{cor}\label{InfinityTerm}

The spectral sequence
$(\tilde{E},\tilde{d})$
has $(p,q)$ term given by
\[\tilde{E}_{1}^{p,q} = H^{p+q-k}(X_{H,0},\mathcal{O}) \otimes H^{k}(X/X_{H,0},\overline{L}).
\]
In fact using Lemmas \ref{someLemma1} and \ref{someLemma2}, we have
\[\tilde{d}_{1} = \tilde{d}_{1} = \cdots = \tilde{d}_{t^{0} - 1} = 0
\]
and so
\[\tilde{E}_{1} = \tilde{E}_{2} = \cdots = \tilde{E_{t^{0}}}.\]
The differential $d_{t^{0}}$ is given (see Lemmas \ref{someLemma1} and \ref{someLemma2})  by \[d_{t^{0}}(\bullet)  =  -\pi l_{t^{0}} \cup \bullet.  \]  Denote by
\[H^{a}(X_{H,0},\mathcal{O})_{l^{0}_{t^0}} \subset H^{a}(X_{H,0},\mathcal{O})
\]
the kernel in  $H^{a}(X_{H,0},\mathcal{O})$ of the map
\[(l^{0}_{t^0} \cup \bullet): H^{a}(X_{H,0},\mathcal{O}) \to H^{a+1}(X_{H,0},\mathcal{O})
.\]

We notice that
\begin {equation}
\label{koseq} H^{a}(X_{H,0},\mathcal{O})_{l^{0}_{t^0}} = \text{Im}\left((l^{0}_{t^0} \cup \bullet): H^{a-1}(X_{H,0},\mathcal{O}) \to H^{a}(X_{H,0},\mathcal{O})\right)
.
\end {equation}

With this notation, the spectral sequence converges and the term $E_{\infty} = E_{t^{0}+1}$  looks like
\def\csg#1{\save[].[ddddddddrrrrr]!C*+<2pc>[F-,]\frm{}\restore}
\[
\xymatrix@=1.5pc{
q=3 & \csg1 \ \ \ \ \ \ \ \vdots & \vdots & \vdots & \vdots & \vdots &\\
q=2 &
H^{2}(X, L)_{l_{t^{0}}}   &
H^{3}(X,L)_{l_{t^{0}}} &
\ddots & H^{t+1}(X,L)_{l_{t^{0}}} &
0 & \cdots \\
q=1 &
H^{1}(X, L)_{l_{t^{0}}}   &
H^{2}(X,L)_{l_{t^{0}}} &
\ddots & H^{t}(X,L)_{l_{t^{0}}} &
0 & \cdots \\
q= 0 & H^{0}(X,L)_{l_{t^{0}}}   &
H^{1}(X,L)_{l_{t^{0}}} &
\ddots & H^{t-1}(X,L)_{l_{t^{0}}} &
0 & \cdots \\
q=-1 &    &
H^{0}(X,L)_{l_{t^{0}}}  &
\ddots & H^{t-2}(X,L)_{l_{t^{0}}} & 0 & \cdots \\
\vdots &
&   &
\ddots & \vdots & 0 & \cdots \\
q=-t^{0}+1 & & & & H^{0}(X,L)_{l_{t^{0}}} & 0 & \ddots \\
q=-t^{0} & & & &  & 0 & \ddots \\
 \vdots & & & &  &  & \ddots \\
 & p= 0 & p= 1 & p= 2 & p=t^{0}-1 & p=t^{0} & \cdots
}
\]
where
\[H^{j}(X,L)_{l_{t^{0}}}  = H^{j-k}(X_{H,0},\mathcal{O})_{l^{0}_{t^{0}}} \otimes H^{k}(X/X_{H,0},\overline{L}) \text{ for } k \leq j \leq g.
\]
Also, observe that the above vector space is zero for $j=k$ and $j>g_{0}+k$.
\end{cor}

\ \hfill $\Box$

\begin{thm}\label{main} Let $\mycal{L}=\mycal{L}_{((H,\chi);l(\hbar)}$ be a line bundle on $\bbX_{\bpi}$.  Assume that $l(\hbar)^{0} \neq 0$ and $\chi|_{X_{H,0}} =1$.

In this case, there is an isomorphism of $\mathbb C[[\hbar]]$ modules
\[H^j(X, \mycal{L}) \cong \mathbb C[\hbar]/(\hbar^{t^{0}}) \otimes_{\mathbb C} (l_{t^{0}} \cup H^{j-1}(X,L)).\]

\end {thm}
\noindent {\bf Proof.}
As the spectral sequence $(E(X, \mycal L), d(X,\mycal L)$ converges, Theorem \ref{specseqthm}, together with the structure Theorem \ref{structure} and Corollary \ref{InfinityTerm} implies that
\[H^j(X, \mycal{L}) \cong \mathbb C[\hbar]/(\hbar^{t^{0}}) \otimes_{\mathbb C} H^{j}(X,L)_{l_{t^{0}}}.\]
By equation (\ref{koseq}), we deduce that
\[H^j(X, L)_{l_{t^{0}}} \cong
({l^{0}_{t^{0}}} \cup H^{j-k-1}(X_{H,0},\mathcal{O})) \otimes H^{k}(X/X_{H,0},\overline{L}).\]
Using equations (\ref{assoc},\ref{CompatOfLs}), the map $s^\vee \otimes \rho^*$ (cf. equation (\ref{identify})) induces an isomorphism
\[({l^{0}_{t^{0}}} \cup H^{j-k-1}(X_{H,0},\mathcal{O})) \otimes H^{k}(X/X_{H,0},\overline{L}) \cong (l_{t^0} \cup H^{j-1}(X,L)),\]

This completes the proof.

\ \hfill $\Box$

\begin{cor}\label{ThreeParts}
For any line bundle $\mycal{L} = \mycal{L}_{((H,\chi);l(\hbar))}$ on $\bbX_{\bpi}$, the following holds:
\begin{enumerate}
\item If $\chi|_{X_{H,0}} \neq 1$ then $H^{j}(X,\mycal{L})= 0$ for all $j$.
\item If $l(\hbar)^{0} =0$ then $H^{j}(X,\mycal{L}) \cong H^{j}(X,\mycal{L}/\hbar)[[\hbar]]$ for all $j$.
\item In all other cases
\[H^j(X, \mycal{L}) \cong \mathbb C[[\hbar]]/(\hbar^{t^{0}}) \otimes_{\mathbb C} (l_{t^{0}} \cup H^{j-1}(X,L)).\]
where
\[
l_{t^{0}} \cup H^{j-1}(X,L) = H^{j-k}(X_{H,0},\mathcal{O})_{l^{0}_{t^{0}}} \otimes H^{k}(X/X_{H,0},\overline{L}).\]
and so
\[\dim_{\mathbb{C}}(H^{j}(X,\mycal{L}))= t^{0}{{g_{0}-1}\choose{j-k-1}}h^{k}(X,L).
\]
\end{enumerate}
\ \hfill $\Box$
\end{cor}
{\bf Proof.}
\begin{enumerate}
\item This is a consequence of the existence of the spectral sequence
\ref{DefOfSpectralSequence} which converges to the cohomology of $\mycal{L}$ and the fact \ref{green} that $H^{j}(X,L)= 0$ for all $j$ and hence $E_{1} = 0$.
\item
This follows from the degeneration (\ref{InfinityTerm}) at $E_{1}$  \ref{DefOfSpectralSequence} and the structure theorem given in Observation \ref{structure}.
\item
This is the content of Theorem \ref{main}.
\end{enumerate}
Although we did not use it, it is interesting to note that the structure of the cohomology of $\mycal{L}/\hbar^{t^{0}}$ and of $\mycal{L}$ are quite different.  Indeed we have the following Lemma.
\begin{lem}\label{descriptions} Let $\mycal{L}=\mycal{L}_{((H,\chi);l(\hbar))}$ be a line bundle on $\bbX_{\bpi}$.  Assume that $l(\hbar)^{0} \neq 0$.  Then
there is an isomorphism of $\mathbb{C}[\hbar]/(\hbar^{t^{0}})$ modules
\begin{equation}\label{someEquation}
H^{j}(X,L)[\hbar]/(\hbar^{t^{0}}) \cong
H^{j}(X,\mycal{L}/(\hbar^{t^{0}}))
\end{equation}
\end{lem}
{\bf Proof.}
The proof is precisely analogous to part (2) of Corollary \ref{ThreeParts}.  There is a spectral sequence converging to the cohomology of $\mycal{L}/(\hbar^{t^{0}})$ as a $\mathbb{C}[\hbar]/(\hbar^{t^{0}})$ module.  An analogue of Corollary \ref{InfinityTerm} shows that all the differentials are zero (because the analogue of $l(\hbar)^{0}$ vanishes for $\mycal{L}/(\hbar^{t^{0}})$) and the term $E_{1}^{p,q} = H^{p+q}(X,L)$, in the range that $Gr(\mycal{L}/(\hbar^{t^{0}})) \neq 0$.  Finally, an appeal to Observation \ref{structure} gives us Equation \ref{someEquation}.
\ \hfill $\Box$

\section{Explicit Cocycle Representatives for Cohomology Classes}
In this section, we find explicit cocycle representatives for a basis of $H^{j}(X,\mycal{L})$.  In order to make this feasible, we restrict to the case that $t^{0} = t$ including the case where they are both infinity.  Let us return to the cases in the Corollary \ref{ThreeParts}.
\begin {enumerate}
    \item
    $\chi|_{X_{H_0}}\neq 1$.
    Then by Corollary \ref{ThreeParts}.1 we have
    $H^{j}(X,\mycal{L}) = 0$ for all $j$ so there is nothing to do here.
    \item
    $l(\hbar)^{0} = 0$ and $t^0=t.$
    In this case we of course have that $t=t^0=\infty$.  Therefore $\mycal{L} = L[[\hbar]]$ and $H^{j}(X,\mycal{L}) \cong H^{j}(X,L)[[\hbar]]$.  Then by Corollary \ref{ThreeParts}.2 together with Lemma \ref{identify} we have
    $$H^{j}(X,\mycal{L}) \cong H^{k}(X/X_{H,0},\overline{L}) \otimes H^{j-k}(X_{H,0},\mathcal{O})[[\hbar]]$$ for all $j$, where $k=k_L$.   Notice that the differential operator $P$ coming from the Poisson structure vanishes on $\rho^{-1} \mathcal{O}_{X/X_{H,0}} \otimes \rho^{-1} \mathcal{O}_{X/X_{H,0}}$ due to equation \ref{condition}.  Explicit representatives are powers of $\hbar$ times the classes in the discussion of Lemma \ref{green} via the map
    \[Z^{k}(\Lambda/(\Lambda \cap V_{H,0}),\mathcal{O}(V/V_{H,0}),\overline{\varphi}) \otimes Z^{i}(\Lambda \cap V_{H,0},\mathcal{O}(V_{H,0}))\otimes \mathbb{C}[[\hbar]] \to Z^{k+i}(\Lambda, \mathcal{O}(V),\Phi).
    \]
    \item In the remainder of this section we look into the case that $\chi|_{X_{H_0}}= 1$ and $l(\hbar)^{0} \neq 0$.

\end {enumerate}
The short exact sequence
    \begin {equation}\label{LScriptLses}
    0 \to \mycal{L} \overset{\hbar^{t^{0}}} \to \mycal{L} \to \mycal{L}/\hbar^{t^{0}} \to 0
    \end {equation}
    gives rise to a long exact sequence in cohomology
    \begin {equation}
    \label {boundmor2}
    \cdots \to H^{j-1}(X,\mycal{L}/\hbar^{t^{0}}) \overset{}  \to H^{j}(X,\mycal{L}) \to H^{j}(X,\mycal{L}) \to H^{j}(X,\mycal{L}/\hbar^{t^{0}}) \to \cdots.
    \end {equation}
\begin{defn}
Let $\mycal{L}$ be a line bundle on $\bbX_{\bpi}$  such that $l(\hbar)^{0} \neq 0.$
We define $\alpha_{\mycal{L}} \in Ext^{1}(X,\mycal{L}/\hbar^{t^{0}},\mycal{L})$ as the extension class of the short exact sequence \ref{LScriptLses}.
\end{defn}
An explicit formula for $\alpha_{\mycal{L}}$ will be
given in formula \ref{explictAlpha}.  With this definition in mind, we have
\begin{lem}\label{cupl2} Let $\mycal{L}$ be a line bundle on $\bbX_{\bpi}$ such that $l(\hbar)^{0} \neq 0$ and $\chi|_{X_{H,0}} =1$.  Then

\[H^{j}(X,\mycal{L}) = \alpha_{\mycal L} \cup H^{j-1}(X,\mycal{L}/\hbar^{t^{0}}).
\]

\end {lem}
{\bf Proof.}
Notice that Theorem \ref{main} shows that $H^{j}(X,\mycal{L})$ is killed by $\hbar^{t^{0}}$ and hence the connecting map
\[H^{j-1}(X,\mycal{L}/\hbar^{t^{0}}) \overset{}  \to H^{j}(X,\mycal{L})
\]
is surjective.  The connecting map is given by the cup product with $\alpha_{\mycal{L}}$.
\ \hfill $\Box$

We have the following commutative diagram
\begin{equation}\label{cupprodcom}
\xymatrix@C+2pc{
Ext^1_{\Ah{X}}(\mycal L/h^{t^{0}}, \mycal L) \times  \ar[d]_-{} H_{\Ah{X}}^{j-1}(X, \mycal L/h^{t^{0}}) \ar[r]^-{\cup}
 &  H_{\Ah{X}}^{j}(X,\mycal L)\ar[d]^-{} \\
Ext^1_{\mathbb{C}_{X}}(\mycal L/h^{t^{0}}, \mycal L) \times  H_{\mathbb{C}_{X}}^{j-1}(X, \mycal L/h^{t^{0}}) \ar[r]_-{\cup} &  H_{\mathbb{C}_{X}}^{j}(X, \mycal L).
}
\end{equation}

Here, the maps on cohomology are isomorphisms.  We notice that for any sheaves of vector spaces $\mycal F, \mycal G$ over $X$, we have
$$Hom_{\mathbb{C}_{X}}(\mycal F,\mycal G) = H^0(\Lambda, Hom_{\mathbb{C}_{V}}(p^{-1}\mycal F, p^{-1} \mycal G)).$$
The Grothendieck spectral sequence for the composition of the functors $H^0(\Lambda,-)$ and \\
$Hom_{\mathbb C}(p^{-1}\mycal F, p^{-1} -)$
is the following convergent first quadrant spectral sequence (cf. \cite{Weibel} Theorem 5.8.3)
\begin {equation}
\label {spse}
E_2^{p,q} := H^p(\Lambda, Ext_{\mathbb{C}_{V}}^q(p^{-1}\mycal F, p^{-1} \mycal G)) \Rightarrow Ext^{p+q}_{\mathbb{C}_{V}}(\mycal F,\mycal G).
\end {equation}

This spectral sequence, applied with $\mycal{F} = \mycal{L}/\hbar^{t^{0}}$ and $\mycal{G} = \mycal{L}$  gives rise to a left exact sequence (cf. \cite{Weibel} Theorem 5.8.3)
\begin{equation} \label{ses1}
\xymatrix @W0pc @C-1pc{
& 0 \ar[r] & H^{1}(\Lambda, Hom_{\mathbb{C}_{V}}(p^{-1}\mycal L/\hbar^{t^{0}},p^{-1}\mycal L)) \ar[r]^-{\beta} & Ext^{1}_{{\mathbb C}_{X}}(\mycal{L}/\hbar^{t^{0}},\mycal{L}) \ar[r]^-{\gamma} & H^{0}(\Lambda, Ext^1_{\mathbb{C}_{V}}(p^{-1}\mycal L/\hbar^{t^{0}},p^{-1}\mycal L)).}
\end {equation}

We want to compute the cup product of $\alpha_{\mycal{L}}$ with elements in $H^{j-1}(X,\mycal L/\hbar^{t^{0}})$ (cf. Lemma \ref{cupl2}).
To this end, we let $\tilde \alpha_{\mycal L}$ denote the image of $\alpha_{\mycal L}$ under the morphism
$$Ext^1_{\Ah{X}}(\mycal L/\hbar^{t^{0}}, \mycal L) \to Ext^1_{\mathbb{C}_{X}}(\mycal L/\hbar^{t^{0}}, \mycal L).$$
Because the sequence
\[0 \to p^{-1}\mycal{L} \to p^{-1} \mycal{L} \to p^{-1} \mycal{L}/\hbar^{t^{0}} \to 0
\]
splits as a sequence of sheaves of vector spaces, we have $\gamma(\tilde\alpha_L)=0$.
Since \ref{ses1} is exact, we deduce that there exists an element
$$a_L \in H^{1}(\Lambda, Hom_{\mathbb{C}_{V}}(\mathcal{O}_{V}[\hbar]/(\hbar^{t^{0}}),
\mathcal{O}_{V}[[\hbar]]),\Psi) $$ such that
$\tilde \alpha_{\mycal L} = \beta(a_{\mycal L})$.
This, together with the commutative diagram \ref{cupprodcom} implies that for $\xi \in H^{j-1}(X,\mycal L/h^{t^{0}})$, we have
$$\alpha_{\mycal L} \cup \xi = a_{\mycal L} \cup \xi$$ as elements in $H^j(X,\mycal L)$.
We now compute $a_{\mycal L}$ in group cohomology
\[H^{1}(\Lambda, Hom_{\mathbb{C}_{V}}(\mathcal{O}_{V}[\hbar]/(\hbar^{t^{0}}),\mathcal{O}_{V}[[\hbar]]),\Psi) \]
where $\Psi$ is the $\Lambda$ action induced on $Hom_{\mathbb{C}_{V}}(\mathcal{O}_{V}[\hbar]/\hbar^{t^{0}},\mathcal{O}_{V}[[\hbar]])$ by the translation actions on $p^{-1}\mycal{L}/\hbar^{t^{0}}$ and $p^{-1}\mycal{L}$.
In order to do this we can consider any \[\mu \in Hom_{\mathbb{C}_{V}}(\mathcal{O}_{V}[\hbar]/(\hbar^{t^{0}}), \mathcal{O}_{V}[[\hbar]])\] such that $\mu$ is the identity modulo $\hbar^{t^{0}}$. Then
\[(a_{\mycal{L}})_{\lambda}(f) = \frac{1}{\hbar^{t^{0}}}(\delta^{\Psi} \mu)_{\lambda}(f).
\]
In order to calculate this let
\begin{equation}\label{regularphi}\phi \in Z^{1}(\Lambda,\mycal{A}_{\bpi}(V)^{\times}/\hbar^{t^{0}})
\end{equation} be the reduction of $\Phi$ modulo $\hbar^{t^{0}}$.
Then for $f \in \mathcal{O}(V)$ we have (note that we will use  ${(A_{\lambda}^{\phi})}^{-1}(f) = (f\phi_{\lambda}) \circ T_{\lambda}^{-1}$)
\begin {equation}
\begin {array}{ll}
(\delta^{\Psi} \mu)_{\lambda} (f) = (A^{\Psi}_{\lambda}(\mu) - \mu)(f)
=A^{\Phi}_{\lambda}(\mu(({{A^{\phi}}_{\lambda}})^{-1} f)) -\mu (f) \cr\cr =  (\mu({(A^{\phi}_{\lambda})}^{-1} f) \circ T_{\lambda}) \star \Phi_{\lambda}^{-1} - \mu(f)
= \mu(f \phi_{\lambda}) \star \Phi_{\lambda}^{-1} - \mu(f)  \cr\cr = \mu(f \phi_{\lambda})\star(\Phi_{\lambda}^{-1} - \phi_{\lambda}^{-1})
\end {array}
\end {equation}
If we chose $\mu$ to be the inclusion with no higher powers of $\hbar$, then we an element \[a_{\mycal{L}} \in Z^{1}(\Lambda, Hom_{\mathbb{C}_{V}}(\mathcal{O}_{V},\mathcal{O}_{V}[[\hbar]]),\Psi)\]
given by
\begin{equation}\label{explictAlpha}
(a_{\mycal{L}})_{\lambda}(f) = (f \phi_{\lambda})\star\frac{(\Phi_{\lambda}^{-1} - \phi_{\lambda}^{-1})}{\hbar^{t^{0}}}.
\end{equation}
Therefore, the product
\[Ext_{\mycal{A}_{X}}^{1}(\mycal{L}/\hbar^{t^{0}},\mycal{L}) \times H^{j-1}(X,\mycal{L}/\hbar^{t^{0}}) \to H^{j}(X,\mycal{L})
\]
evaluated on the pair $(\alpha_{\mycal{L}},\xi)$ takes on the form
\[(a_{\mycal{L}}\cup \xi)_{\lambda_{0},\dots,\lambda_{j}} =  \frac{(\delta^{\Psi}\mu)_{\lambda_{0}} }{\hbar^{t^{0}}}(A_{\lambda_{0}}^{\phi}(\xi_{\lambda_{1},\dots,\lambda_{j}} )) =
\frac{(\delta^{\Psi} \mu)_{\lambda_{0}}}{\hbar^{t^{0}}}({\phi_{\lambda_{0}}}^{-1}(\xi_{\lambda_{1},\dots,\lambda_{j}} \circ T_{\lambda_{0}}))
\]
\[= ((\phi_{\lambda_{0}}\phi_{\lambda_{0}}^{-1})(\xi_{\lambda_{1},\dots,\lambda_{j}} \circ T_{\lambda_{0}})) \star \frac{(\Phi_{\lambda_{0}}^{-1} - \phi_{\lambda_{0}}^{-1})}{\hbar^{t^{0}}}.
\]
We deduce
\begin{equation} \label{AlphaCupEqn}
(a_{\mycal{L}}\cup \xi)_{\lambda_{0},\dots,\lambda_{j}} = (\xi_{\lambda_{1},\dots,\lambda_{j}} \circ T_{\lambda_{0}}) \star \frac{(\Phi_{\lambda_{0}}^{-1} - \phi_{\lambda_{0}}^{-1})}{\hbar^{t^{0}}}.
\end{equation}

Notice one can see independently that $\alpha_{\mycal{L}}\cup \xi \in Z^j(\Lambda, \mathcal O(V)[[\hbar]], \Phi):$
By analogy to equation \ref{form} we deduce that
$$\alpha_{\mycal{L}}\cup \xi = \frac {\delta^\Phi(\xi)} {\hbar^{t^{0}}}.$$
This implies that
$$\delta^\Phi(\alpha_{\mycal{L}}\cup \xi) = \frac {(\delta^\Phi)^2(\xi)} {\hbar^{t^{0}}} = 0.$$

\begin{observation}\label{AlphaObservation}
The image of $\alpha_{\mycal L}$ under the modulo $\hbar$ reduction
     $$Ext^{1}(\mycal{L}/\hbar^{t^{0}},\mycal{L}) \to Ext^{1}(L,L) = H^1(X,\mathcal O)$$ is equal to $-\pi l_{t^{0}}$.
\end{observation}

\begin{observation}
 Let $\mycal{L}=\mycal{L}_{((H,\chi);l(\hbar))}$ be a line bundle on $\bbX_{\bpi}$.  Assume that $l(\hbar)^{0} \neq 0$ , $t^{0} = t$, and that $\chi|_{X_{H,0}} =1$.  In this case we simply have $\mycal{L}/\hbar^{t^{0}} = L[\hbar]/\hbar^{t^{0}}$.
One can check using
\[\bpi(d (\rho^{-1} \mathcal{O}_{X/X_{H,0}}) , d (\rho^{-1} \mathcal{O}_{X/X_{H,0}})) = 0
\]
which follows from Equation (\ref{condition}),
 that the isomorphism in Lemma \ref{descriptions} is induced by the map
\[Z^{j}(\Lambda,\mathcal{O}(V),\varphi)[\hbar]/(\hbar^{t^{0}}) \to Z^{j}(\Lambda,\mathcal{O}(V)[\hbar]/(\hbar^{t^{0}}),\phi)
\]
given by the identity on $C^{j}(\Lambda,\mathcal{O}(V)[\hbar]/\hbar^{t^{0}})$.
In other words, the above is well defined and gives an isomorphism in cohomology because the $\star$-products involved agree with commutative products for functions in $\mathcal{O}(V/V_{H,0})$ with the action $A^{\overline{\varphi}}$ where \[\varphi_{\lambda} = \overline{\varphi}_{\rho(\lambda)} \circ \rho.\]
\end{observation}
\begin{lem}\label{AlphaIsom} Let $\mycal{L}=\mycal{L}_{((H,\chi);l(\hbar))}$ be a line bundle on $\bbX_{\bpi}$.  Assume that $l(\hbar)^{0} \neq 0$ and that $\chi|_{X_{H,0}} =1$.
Then we have an isomorphism of $\mathbb{C}[[\hbar]]$ modules given by $\alpha_{\mycal{L}} \cup \bullet$
\[
(\mathbb{C}[\hbar]/\hbar^{t^{0}}) \otimes  H^{k}(X/X_{H,0},\overline{L}) \otimes (H^{j-k-1}(X_{H,0},\mathcal{O})/<l^{0}_{t^{0}}>) \to H^{j}(X,\mycal{L})
\]
where we implicitly used the identification in \ref{identify} in order to apply $\alpha_{\mycal{L}}$ to $ H^{k}(X/X_{H,0},\overline{L}) \otimes H^{j-k-1}(X_{H,0},\mathcal{O})$.  In particular, the cohomology is non-zero only the range $k+1 \leq j \leq k+g_{0}$.
\end{lem}

\noindent {\bf Proof.}
\
Since we know from \ref{ThreeParts} that $H^{k}(X,\mycal{L})$ is $t^{0}$-torsion, we see that $H^{j-1}(X,\mycal{L}/(\hbar^{t^{0}}))$ surjects onto $H^{j}(X,\mycal{L})$ via the connecting map (cup product with the extension class),
\[H^{j}(X,\mycal{L}) = \alpha_{\mycal{L}} \cup H^{j-1}(X,\mycal{L}/(\hbar^{t^{0}})) \cong \alpha_{\mycal{L}} \cup  H^{j-1}(X,L) \otimes \mathbb{C}[\hbar]/(\hbar^{t^{0}})
\]
The calculation of the kernel of the map given by the cup product with $\alpha_{\mycal{L}}$ follows from Observation \ref{AlphaObservation}.
\
\ \hfill $\Box$
\begin{cor}  Let $\mycal{L}=\mycal{L}_{((H,\chi);l(\hbar))}$ be a line bundle on $\bbX_{\bpi}$.  Assume that $l(\hbar)^{0} \neq 0$ , $t^{0} = t$, and that $\chi|_{X_{H,0}} =1$.  We take $j$ such that $k+1 \leq j \leq k+g_{0}$ because this is the only range in which there is non-zero cohomology.
Let $a^{1}, \dots, a^{g_{0}-1}$ be a basis of a compliment to $<l^{0}_{t^{0}}>$ in $\overline{V_{H,0}}^{\vee}$.  Let $\{b^{r}\}$ be elements of $Z^{k}(\Lambda/(\Lambda \cap V_{H,0}),\mathcal{O}(V/V_{H,0}),\overline{\varphi})$ whose cohomology classes are a basis for \[H^{k}(\Lambda/(\Lambda \cap V_{H,0}),\mathcal{O}(V/V_{H,0}),\overline{\varphi}) \cong H^{k}(X/X_{H,0},\overline{L})\] then if we define
\[a^{I} = a^{i_{1}} \cup \cdots \cup a^{i_{j-k-1}}
\]
the following
\begin{equation}\label{finalExplicit}
\hbar^{c}\left(b^{r}_{\rho(\lambda_{1}), \dots , \rho(\lambda_{k})} \circ \rho) \circ T_{\lambda_{0}}\right) \left(s^{\vee}a^{I}_{\lambda_{k+1}, \dots , \lambda_{j-1}}\right)  \phi_{\lambda_{0}}^{-1}\frac{1}{\hbar^{t^{0}}}\left(\exp \left(-\sum_{m=t^{0}}^{\infty} \hbar^{m}\pi l_{m} (\lambda_0) \right)-1\right)
\end{equation}
for $0 \leq c < t^{0}$ , $1 \leq r \leq h^{k}(X,L)$ and $1 \leq i_{1} < \cdots < i_{j-k-1} \leq g_{0}-1$ are  elements in $Z^{j}(\Lambda, \mathcal{O}(V)[[\hbar]],\Phi)$ whose cohomology classes are a basis for \[H^{j}(\Lambda, \mathcal{O}(V)[[\hbar]],\Phi) \cong H^{j}(X,\mycal{L}).\]

\end{cor}

{\bf Proof.}
Notice that we can rewrite
\[\Phi_{\lambda}^{-1} = \varphi_{\lambda}^{-1} \exp\left(- \pi \sum_{m=1}^{\infty}\hbar^{m}l_{m}(\lambda) \right) = \left(\overline{\varphi}_{\rho(\lambda)}^{-1} \circ \rho \right) \exp \left(- \pi \sum_{m=1}^{\infty}\hbar^{i}l_{m}(\lambda) \right)
\]
and similarly
\[\phi_{\lambda}^{-1} = \varphi_{\lambda}^{-1}\exp \left(- \pi \sum_{m=1}^{t^{0} - 1}\hbar^{m}l_{m}(\lambda)\right)=\left(\overline{\varphi}_{\rho(\lambda)}^{-1} \circ \rho \right) \exp\left(- \pi \sum_{m=1}^{t^{0} - 1}\hbar^{m}l_{m}(\lambda)\right).
\]
Since by equation \ref{condition} the bi-differential operator $P$ vanishes on $\rho^{-1}\mathcal{O}_{X/X_{H,0}} \otimes \rho^{-1}\mathcal{O}_{X/X_{H,0}}$ we have
\[\left((b^{i}_{\rho(\lambda_{1}),\dots, \rho(\lambda_{k})} \circ T_{\rho(\lambda_{0})}) \circ \rho \right)  \star  \left( \overline{\varphi}_{\rho(\lambda)}^{-1} \circ \rho \right) = \left( (b^{i}_{\rho(\lambda_{1}),\dots, \rho(\lambda_{k})} \circ T_{\rho(\lambda_{0})}) \circ \rho  \right) \left( \overline{\varphi}_{\rho(\lambda)}^{-1} \circ \rho \right).
\]
Therefore, the existence of a collection as in \ref{finalExplicit} follows from Corollary \ref{AlphaIsom} and Equation \ref{AlphaCupEqn}.
\ \hfill $\Box$
\begin{rem}\label{exception}
In the case $t<t^{0}< \infty$ one can still identify cocycles representing a partial basis of $H^{j}(X,\mycal{L})$ using the cup product with the extension class of
\[0 \to \mycal{L} \overset{\hbar^{t}} \to \mycal{L} \to \mycal{L}/\hbar^{t} \to 0.
\]
\end{rem}
\section{Appendix 1}
In order to do this computation we will use the canonical isomorphism \[C: p^{-1}\mycal{L} \cong \mathcal{O}_{V}[[\hbar]] \] (to be explained below) to relate the translation action on $p^{-1}\mycal{L}(V)$ to some action which we compute on $\mathcal{O}_{V}(V)[[\hbar]]$.

We have a canonical trivialization of the pullback of $\mycal{L}$ given by composing the pullback of the inclusion $\mycal{L} \in p_{*}\mathcal{O}_{V}[[\hbar]]$ with the canonical morphism $p^{-1}p_{*} \mathcal{O}[[\hbar]] \to \mathcal{O}[[\hbar]]$.
\[p^{-1} \mycal{L} \subset p^{-1}p_{*} \mathcal{O}[[\hbar]] \to \mathcal{O}[[\hbar]].
\]
Let us call the combined isomorphism \[C:p^{-1} \mycal{L} \to \mathcal{O}_{V}[[\hbar]].\]  In particular, we use the same letter for the map on global sections \[C: H^{0}(V,p^{-1} \mycal{L}) \to H^{0}(V, \mathcal{O}_{V}[[\hbar]]).\]

In order to derive a specific formula for $C$, consider the element $F \in H^{0}(V, p^{-1} p_{*} \mathcal{O}_{V}[[\hbar]])$ defined by
\[F(v)(w) = \Phi_{w-v}(v) .
\]
We claim in fact that $F$ is a nowhere vanishing global section of $p^{-1} \mycal{L}$.  In order to show this, it suffices to prove that $F$ locally belongs to $p^{-1} \mycal{L}$.  In other words we want to see that for every $\lambda \in \Lambda$ that,
\[F(v)(w) \star \Phi_{\lambda}(w) = F(v)(w) \circ t_{\lambda}
\]
for $p(v) = p(w)$ lying in $U$.

The crucial point here, is that $T$ acts on the $v$ variable {\it horizontally} and $t$ acts on the $w$ variable {\it vertically}.
\[F(v)(w) \star \Phi_{\lambda}(w) = \Phi_{w-v}(v) \star \Phi_{\lambda}(w) = \Phi_{w-v}(v) \star (\Phi_{\lambda}(v) \circ T_{w-v}) = \Phi_{w-v+\lambda}(v)
\]
and
\[F(v)(w) \circ t_{\lambda} = \Phi_{w-v}(v) \circ t_{\lambda} = \Phi_{w-v+\lambda}(v).
\]

The maps $T_{\lambda *}p^{-1} p_{*} \mathcal{O} \to p^{-1} p_{*} \mathcal{O}$ and $T_{\lambda *}p^{-1} \mycal{L} \to p^{-1} \mycal{L}$ induce the actions on $H^{0}(V,p^{-1} p_{*} \mathcal{O})$ and $H^{0}(V,p^{-1} \mycal{L})$ given by
\[f \mapsto f \circ T_{\lambda}.
\]
 We claim that $C^{-1}(g) = g \star F$  and $C(h) = h \star F^{-1}$.
Indeed, for $w=v+\lambda$
\[C^{-1}(1) = \Phi_{\lambda} \circ t_{-\lambda} = \Phi_{w-v}(w-\lambda) = \Phi_{w-v}(v) = F(v)(w)
\]
We will use $C$ to transport the translation action (defined by $A_{\lambda}(f) = f \circ T_{\lambda}$)
\[A_{\lambda}: T_{\lambda *} p^{-1} \mycal{L} \to p^{-1} \mycal{L}
\]
satisfying
$A_{\lambda_{1}} \circ (T_{\lambda_{1} *}A_{\lambda_{2}}) = A_{\lambda_{1}+\lambda_{2}}$
to the action
\[  A^{\Phi}_{\lambda} = C \circ A_{\lambda} \circ T_{\lambda *} C^{-1} : T_{\lambda *} \mathcal{O}_{V}[[\hbar]] \to \mathcal{O}_{V}[[\hbar]].
\]

\begin{equation}\label{BigAlittleA}
\xymatrix@C+2pc{
T_{\lambda *}  \mathcal{O}_{V} \ar[d]_-{A^{\Phi}_{\lambda}}
 & T_{\lambda *} p^{-1} \mycal{L} \ar[d]^-{A_{\lambda}} \ar[l]^-{T_{\lambda *}C}\\
\mathcal{O}_{V} \ar[r]_-{C^{-1}} &  p^{-1} \mycal{L}.
}
\end{equation}

Let $g = g(v)$ be in $H^{0}(V,\mathcal{O}_{V}[[\hbar]])$.  We compute
\[((g \star F) \circ T_{\lambda}) \star F^{-1} = (g \circ T_{\lambda}) \star (F \circ T_{\lambda}) \star F^{-1} = (g \circ T_{\lambda}) \star (\Phi_{w-v-\lambda} \circ T_{\lambda})  \star \Phi_{w-v}^{-1} = (g \circ T_{\lambda}) \star \Phi_{\lambda}^{-1}.
\]
Therefore
\begin{equation}A^{\Phi}_{\lambda}(g) = (g \circ T_{\lambda}) \star \Phi_{\lambda}^{-1}
.\end{equation}
In other words the action $A^{\Phi}_{\lambda}$ on $H^{0}(V, \mathcal{O}_{V} [[\hbar]])$ is given by
\begin {equation}
\label {trac2}
g(v) \mapsto g(v +\lambda) \star \Phi_{\lambda}(v)^{-1}.
\end {equation}
We can check directly that this is indeed an action:
\[(A^{\Phi}_{\lambda_{1}} A^{\Phi}_{\lambda_{2}} (g)) (v) = g(v + \lambda_{2} + \lambda_{1}) \star \Phi_{\lambda_{2}}(v + \lambda_{1})^{-1} \star \Phi_{\lambda_{1}}(v)^{-1} \]
\[= g(v + \lambda_{2} + \lambda_{1}) \star \Phi_{\lambda_{1}+\lambda_{2}}(v )^{-1} = (A^{\Phi}_{\lambda_{1} + \lambda_{2}} (g))(v).
\]

For a small open set $U \subset V$ and $g \in \mathcal{O}(U)[[\hbar]]$ note that \[C^{-1}(g) \in p^{-1}\mycal{L}(U) \subset \mathcal{O}(p^{-1}(p(U)))[[\hbar]]\] is determined uniquely by the property $C^{-1}(g) \circ t_{\lambda} = C^{-1}(g)  \star \Phi_{\lambda} $.  In other words, for $W \subset p^{-1}(p(U))$ such that $W=\lambda U$ we see that
\[C^{-1}(g)|_{W} = (  g \star \Phi_{\lambda}) \circ t_{-\lambda}.
\]
and
\[C(h_{W}) = (h_{W} \circ t_{\lambda}) \star \Phi_{\lambda}^{-1}.
\]

Pushing the diagram \ref{BigAlittleA} forward and using the identification $p_{*}T_{\lambda *} = p_{*}$ we get
\[
\xymatrix@C+2pc{
p_{ *}  \mathcal{O}_{V} \ar[d]_-{p_{*}A^{\Phi}_{\lambda}}
 & p_{*} p^{-1} \mycal{L} \ar[d]^-{p_{*}A_{\lambda}} \ar[l]^-{p_{ *}C}\\
p_{*}\mathcal{O}_{V} \ar[r]_-{p_{*}C^{-1}} &   p_{*}p^{-1} \mycal{L}
}
\]
This identifies $\mycal{L}$, the $\Lambda$ invariants of $p_{*} p^{-1} \mycal{L}$ with the translation action on $V$ with the $\Lambda$ invariants of $p_{*}\mathcal{O}_{V}$ with the $A^{\Phi}$ action, which in fact agrees with our previous description of $\mycal{L}$ given in equation \ref{eq:DefOfL}.

\begin{rem}\label{ARemark}  Similarly, if $\Phi_{s}$ is the reduction of $\Phi$ modulo $\hbar^{s}$ we denote by $A^{\Phi_{s}}$ the induced action \begin{equation}A^{\Phi_{s}}_{\lambda}(g) = (g \circ T_{\lambda}) \star {(\Phi_{s})}_{\lambda}^{-1}
.\end{equation}on $\mathcal{O}(V)[\hbar]/\hbar^{s}$ and we can replace $H^{i}(X,\mycal{L}/\hbar^{s})$ by
\[H^{i}(\Lambda, \mathcal{O}(V)[\hbar]/\hbar^{s}, \Phi_{s})\] the cohomology groups of $C^{i}(\Lambda, \mathcal{O}(V)[\hbar]/\hbar^{s})$ with the differential $\delta^{\Phi_{s}}$. In particular, for $s=1$ we have $\Phi_{s} = \varphi$ and for $s=t^{0}$ we have $\Phi_{s} = \phi$.

\end{rem}
\section {Appendix 2}
A (decreasing) \emph{filtration} $F$ on a chain complex $(C,d)$ is an ordered family of chain subcomplexes
$$\cdots \subseteq F_{p+1} C \subseteq F_p C \subseteq F_{p-1} \cdots$$
of $C$. For such a filtration there is an associated spectral sequence.
This appears in detail in \cite{Weibel} section 5.4.
In this appendix we give explicit formulae for $E_r, d_r$.
This is taken from loc. cit. 5.4.6.
We omit the bookkeeping subscript $q$, and write $\eta_p$ for the surjection $F_p C \to F_p C/F_{p+1} C = E_0^p$. Next, Weibel introduces
$$A_p^r = \{c \in F_p C: d(c) \in F_{p+r}C\},$$
the elements of $F_p C$ that are cycles modulo $F_{p+r} C$ and their images $Z_p^r = \eta_p(A_p^r)$ in
$E_p^0$ and $B_{p+r}^{r+1} = \eta_{p+r}(d({A_p}^r) )$ in $E_{p+r}^0$. Using this indexing the $Z_p^r$ and $B_p^r = \eta_p(d(A_{p-r+1}^{r-1}) )$ are subobjects of $E_p^0$.
The following holds:
\begin {enumerate}
    \item
    $A_p^r \cap F_{p+1} C = A_{p+1}^{r-1}$.
    \item
    $Z_p^r \cong A_p^r/A_{p+1}^{r-1}$.
    \item
    Hence
    $$E_p^r = \frac {Z_p^r} {B_p^r} \cong \frac {A_p^r+F_{p+1}(C)} {d(A_{p-r+1}^{r-1}) + F_{p+1}(C)} \cong \frac {A_p^r} {d(A_{p-r+1}^{r-1}) + A_{p+1}^{r-1}}.$$
    \item
    $d_p^r: E_p^r \to E_{p+r}^r$ is induced by the differential $d$.
    That is, given $\xi \in E_p^r \cong \frac {A_p^r} {d(A_{p+r-1}^{r-1}) + A_{p+1}^{r-1}}$ we first lift it to any $\Xi \in A_p^r$ whose class in $E_p^r$ is $\xi$. Next, compute $d(\Xi) \in  F_{p+r} C$.
    As $d^2(\Xi)=0$, it follows that $d(\Xi)$ in fact belongs to $A_{p+r}^r$. Finally, its image under $\eta_{p+r}$ gives $d_p^r(\xi) \in E_{p+r}^r$.
\end{enumerate}



\newcommand{\etalchar}[1]{$^{#1}$}

\

\bigskip

\noindent
Department of Mathematics \\
University of Haifa \\
Mount Carmel, Haifa, 31905, Israel \\
\\
\noindent
Landau Center for Mathematical Analysis \\
Einstein Institute of Mathematics \\
Edmond J. Safra Campus, Givat Ram \\
The Hebrew University of Jerusalem \\
Jerusalem, 91904, Israel

\

\bigskip

\noindent
email: ben-bassat@math.haifa.ac.il, noamso@math.huji.ac.il

\end{document}